\documentclass{article}
\usepackage{emsnewsletter}
\usepackage{eurosym}
\newtheorem*{corollary*}{\textbf{Corollary}}
\newtheorem{cor}{\textbf{Corollary}}

\def\bitcoinA{%
  \leavevmode
  \vtop{\offinterlineskip 
    \setbox0=\hbox{B}%
    \setbox2=\hbox to\wd0{\hfil\hskip-.03em
    \vrule height .3ex width .15ex\hskip .08em
    \vrule height .3ex width .15ex\hfil}
    \vbox{\copy2\box0}\box2}}

\author{Cyril Grunspan (De Vinci Research Center, Paris, France)\\  Ricardo P\'erez-Marco (CNRS, IMJ-PRG, Sorbonne Universit\'e, Paris, France)}
\title{The mathematics of Bitcoin}

\begin{document}
   
\maketitle

\begin{abstract}

\end{abstract}


\section{Introduction to Bitcoin.}

Bitcoin is a new decentralized payment network that started operating in January 2009. This new technology 
was created by a pseudonymous author, or group of authors, called Satoshi Nakamoto in an article that was 
publically released \cite{Nakamoto08} in the cypherpunk mailing list. The cypherpunks are anarchists 
and cryptographers that who have been concerned with personal privacy in the Internet  since the 90's.
This article follows on a general presentation of Bitcoin by the second author \cite{PerezMarco2016}. 
We refer to this previous article for general background. 
Here we focuss on mathematics being a feature of the security and effectiveness of Bitcoin protocol.

Roughly speaking the Bitcoin's protocol is a mathematical algorithm on a network which manages transaction data 
and builds majority consensus among the participants. Thus, if a majority of 
the participants are honest, then we get an honest \textit{automatic} 
consensus. Its main feature is \textit{decentralization}, which means that 
no organization or central structure is in charge. 
The nodes of the network are voluntary participants that enjoy equal rights and obligations. The network is open 
and anyone can participate. Once launched the network is resilient and unstopable. 
It has been functioning permanently without significant interruption since january 2009.

The code and development are open. The same code has been reused and modified to create hundreds of other cryptocurrencies based on the same principles. The security of the network relies on strong cryptography 
(several orders of magnitude stronger than the cryptography used in classical financial services). 
For example, classical hash functions (SHA256, RIPEMD-160) and elliptic curve digital signatures algorithm 
(ECDSA) are employed. The cryptography used is very standard and well known, so we will dwell
on the mathematics of these cryptographic tools, but interesting cryptographical research is 
motivated by the special features of other cryptocurrencies.

\medskip
\textbf{Nodes and mining.}
\medskip

The Bitcoin network is composed by nodes, that correspond to the Bitcoin program running on different machines, 
that communicate with their neighbourds. Properly formatted Bitcoin transactions flood the network, and 
are checked, broadcasted and validated continuously by the nodes which follow a precise set of rules. There is no way 
to force the nodes to follow these rules. Incentives are created 
so that any divergence from the rules is economically penalised, thus creating a virtuous cycle. In this way, 
the network is a  complex Dynamical System and it is far from obvious that it is stable. The 
stability of this system is a very interesting and fundamental mathematical problem. In its study we will 
encounter special functions, martingale theory, Markov chains, Dyck words, etc.

\begin{figure}
\centering
\includegraphics[width=.6\linewidth]{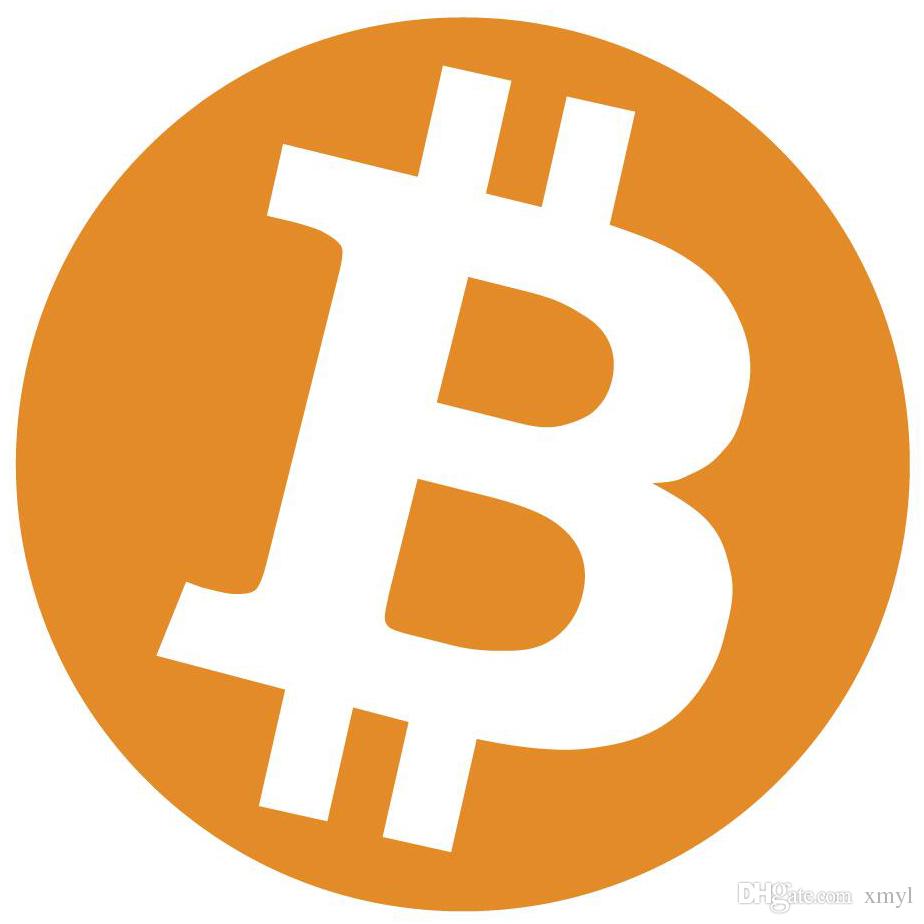}

\medskip

\caption{The Bitcoin logo.}\label{fig:bitcoin}

	
\end{figure}

Nodes in the network broadcast transactions and can participate in their validation. The process of 
validating transactions is also called ``mining'' because it is related to the production of new bitcoins. 
The intuition behind Bitcoin is that of a sort of ``electronic gold'' and the rate of production of bitcoins 
is implemented in the protocol rules. On average every 10 minutes a block of transactions is validated and new bitcoins 
are minted in a special transaction without bitcoin input, called the coinbase transaction. 
At the beginning $50 \, \bitcoinA$ were created in each block, and about every 4 years (more precisely, every $210\, 000$
blocks), the production is divided by $2$. This event is called a ``halving''. So far, we had two halvings, and 
the production is currently set at $12.5\, \bitcoinA$  per $10$ minutes, or $1\, 800\, \bitcoinA$ per day. The next halving will occur 
on May 2020. This geometric decrease of the production limits the total amount of bitcoins to $21$ millions. 
Currently, about $18$ millions have already been created. Each block containing the validated transactions can contain about 3 to 4
thousand transactions and has a size of about 1 Mb. These blocks are linked together cryptographically, and the set 
of all these blocks forms the ``blockchain'' that contains the full history of Bitcoin transactions. This data is 
stored efficiently, and the current blockchain is only about $260.000$ Mb. The cryptographical link between blocks 
is provided by the mining/validation procedure that is based on a hash function and a ``Proof-of-Work''. It costs 
computation power to validate a block and this is what ensures that the data cannot be tampered or corrupted. In 
order to modify a single bit of a block, we must redo all computations that has been used to build all the 
subsequent blocks until the last current one. Currently the computation power needed to change only the last 
few blocks of the more than $600$ thousands composing the  
blockchain is beyond the capabilities of any state or company.

The mining/validation procedure is a sort of decentralized lottery. A miner (this is a node engaging in validating
transactions) packs together a block of floating not yet validated transactions, and builds a header of this block
that contains a hash of the previous block header. The hash algorithm used is SHA-256 (iterated twice) that outputs 
256 bits. Mathematically, a hash function is a deterministic one way function: It is easy to compute, but practically
impossible to find pre-images, or collisions (two files giving the same output). Also it enjoys pseudo-random 
properties, that is, if we change a bit of the input, the bits of the output behave as uncorrelated random variables 
taking the values 0 and 1 with equal probabilities. The mining procedure consists of finding a hash that is below 
a pre-fixed threshold, which is called the \textit{difficulty}. The difficulty is updated every two weeks (or more 
precisely every 2016 blocks) so that the rate of validation remains at 1 block per 10 minutes. The pseudo-random 
properties of the hash function ensure that the only way to find this hash is to 
probe many hashes by changing a parameter in the header (the nonce). The first miner to find a solution makes 
the block public, and the network 
adopts the block as the last block in the blockchain. 

It can happen that two blocks are simultaneously validated in different parts of the network. Then a competition 
follows between the two candidates, and the first one to have a mined block on top of it wins. The other one is discarded and is called an \textit{orphan} block. The blockchain with 
the larger amount of work (which is in general the longer one) is adopted by the nodes.

When a transaction is included in the last block of the blockchain, we say that it has one confirmation. Any extra block 
mined on top of this one gives another confirmation to the transaction and engraves it further inside the blockchain.

This apparently complicated procedure is necessary to ensure that neither the network nor the blockchain cannot be corrupted. 
Any participant must commit some computer power in order to participate in the decision of validation. 
The main obstacle for the invention of a decentralised currency was 
to prevent double spend without a central accounting authority. 
Hence, the first mathematical problem than Nakamoto considers in his founding article \cite{Nakamoto08} is to estimate the 
probability of a double spend. In the following we consider this and other stability problems, and prove 
mathematically the (almost general) stability of the mining rules.

\section{The mining model.}

We consider a miner with a fraction $0<p\leq 1$ of the total hashrate. His profit comes from the block rewards of his validated blocks. 
It is important to know the probability of success when mining a block. The average number of blocks per unit of time 
that he succeeds mining is proportional to his 
hashrate $p$. The whole network takes on average $\tau_0=10 \text{ min}$ to validate a block, hence our miner takes on average $t_0=\frac{\tau_0}{p}$.
We consider the random variable $\mathbf{T}$ giving the time between blocks mined by our miner. The pseudo-random properties of the hash function
shows that mining is a Markov process, that is, memoryless. It is then an elementary exercise to show from this property that $\mathbf{T}$ follows 
an exponential distribution,
$$
f_{\boldsymbol{\tau}} (t) =\alpha e^{-\alpha t} \,
$$
where $\alpha =1/t_0 = 1/\mathbb{E} [ \mathbf{T}]$. If the miner starts mining at $t=0$, and if we denote $\mathbf T_1$ the time needed to mine a first block,
then $\mathbf T_2, \ldots ,\mathbf T_n$ the inter-block mining times of successive blocks, then the Markov property shows that the random variables
$\mathbf T_1, \mathbf T_2, \ldots ,\mathbf T_n$ are independent and are all identically distributed following the same exponential law. Therefore, 
the time needed to discover $n$ blocks is 
$$
\mathbf S_n = \mathbf T_1+\mathbf T_2+\ldots +\mathbf T_n \ .
$$
The random variable $\mathbf S_n$ follows the $n$-convolution of the exponential distribution 
and, as is well known, this gives a Gamma distribution with 
parameters $(n,\alpha)$,
$$
f_{\mathbf S_n} (t)= \frac{\alpha^n}{(n-1)!} t^{n-1}e^{-\alpha t} 
$$
with cumulative distribution 
$$
F_{\mathbf S_n}(t)=\int_0^t f_{\mathbf S_n}(u) du= 1-e^{-\alpha t} \sum_{k=0}^{n-1} \frac{(\alpha t)^k}{k!} \ .
$$
From this we conclude that if $\mathbf N(t)$ is the process counting the number of blocks validated at time $t>0$, 
$\mathbf N(t)=\max\{n\geq 0; \mathbf S_n <t\}$, then we have
$$
\mathbb P[\mathbf N(t)=n]=F_{\mathbf S_n}(t)-F_{\mathbf S_{n+1}}(t)=\frac{(\alpha t)^n}{n!}\ e^{-\alpha t} \  ,
$$
and $\mathbf N(t)$ follows a Poisson law with mean value $\alpha t$. This result is classical, and the mathematics 
of Bitcoin mining, as well as other crypto-currencies with validation based on proof of work, are Poisson distributions.

\section{The double spend problem.} \label{sec:double_spend}

The first crucial mathematical problem that deserves attention in the Bitcoin protocol
is the possibility of realisation of a double spend. This was the major obstacle to overcome for the invention of decentralized cryptocurrencies, 
thus it is not surprising that Nakamoto addresses this problem in Section 11 of his founding article \cite{Nakamoto08}. He considers the 
situation where a malicious miner makes a payment, then in secret tries to validate a second conflicting 
transaction in a new block , from the same address, but to a new address that he controls, which allows him to recover the funds.

For this, once the first transaction has been validated in a block in the official blockchain and the vendor delivered the goods (the vendor will not deliver 
unless some confirmations are visible), the only possibility consists 
in rewriting the blockchain from that block. This is feasible if he controls a majority of the hashrate, that is, if his relative hashrate $q$ satisfies $q>1/2$, because then 
he is able to mine faster than the rest of the network, and he can rewrite the last end of the blockchain as he desires. This is the reason why 
a decentralised mining is necessary so that no one controls more than half of the mining power. But even when $0<q<1/2$ he can try to attempt a double 
spend, and will succeed with a non-zero probability. The first mathematical problem consists of computing the probability that the malicious 
miner succeeds in rewriting the last  $n\geq 1$ blocks. We assume that the remaining relative hashrate, $p=1-q$, consists of honest miners 
following the protocol rules.

This problem is similar to the classical gambler's ruin problem. Nakamoto observes that the probaility to catch-up $n$ blocks is 
$$
q_n=\left (\frac{q}{p}\right )^n  \ \ \ \text{(Nakamoto)}
$$
The modelization of mining shows that the processes $\mathbf N(t)$ and $\mathbf N'(t)$ counting the number of mined blocks at time $t$ 
by the honest and malicious miners, respectively, are independent Poisson processes with respective parameters $\alpha$ et $\alpha'$
satisfying
\begin{equation*}
 p=\frac{\alpha}{\alpha+\alpha'} \ ,\ \ \ \  q=\frac{\alpha'}{\alpha+\alpha'}\ .
\end{equation*}
The random variable $\mathbf X_n=\mathbf N'(\mathbf S_n)$ of the number of blocks mined by the attacker when the honest miners have mined their 
$n$-th block follows a negative binomial variable with parameters $(n,p)$ (\cite{GPM1}), thus, for an integer $k\geq 1$ we have 
$$
\mathbb P[ \mathbf X_n = k] = p^kq^k\binom{k+n-1}{k}  \ .
$$
Nakamoto, in section 11 of \cite{Nakamoto08}, approximates abusively $\mathbf X_n=\mathbf N'(\mathbf S_n)$ by $\mathbf N'(t_n)$ where
$$
t_n =\mathbb E[\mathbf S_n]=n\mathbb E[\mathbf T] = \frac{n\tau_0}{p}
$$
This means that he considers the classical approximation of a negative binomial variable by a Poisson variable. Rosenfeld observes in \cite{R14}
that the negative binomial variable seems to be a better approximation. We proved in \cite{GPM1} that this was indeed the case  and we could 
find the exact formula for the double spend probability after $z$ confirmations ($z$ is the classical notation used by Nakamoto for the 
number of confirmations).

\begin{theorem} [\cite{GPM1}, 2017]\label{thm:probability}
After $z$ confirmations, the probability of success of a double spend by attackers with a relative hashrate of $0<q<1/2$ is
 $$
 P(z)=I_{4pq}(z,1/2) \ 
 $$
where $I_x(a,b)$ is the incomplete regularized beta  function
 $$
 I_x(a,b)=\frac{\Gamma(a+b)}{\Gamma(a) \Gamma(b)} \int_0^x t^{a-1}(1-t)^{b-1} \, dt \ .
 $$
\end{theorem}

Bitcoin security depends on this probability computation. It is not just a theoretical result. It allows the estimation of the 
risk of a transaction to be reversed and the number of confirmations required to consider it definitive. For example, if $q=0.1$, after $6$
confirmations, the probability of a double spend is smaller than 1\% (for complete tables see \cite{GPM2}).

In his founding article, Nakamoto tries to compute this probability from his approximate argument and runs a numerical simulation. 
He convinced himself that the probability converges exponentially to $0$ when the number 
of confirmations $z$ goes to infinite (as he states ``we can see the probability drop off exponentially with $z$''). The numerical 
simulation is not a proof, but we can read this statement repeated over and over, but never proved 
before 2017. With the previous exact formula, using classical methods (Watson Lemma), we can prove the following Corollary:

\begin{cor} Let $s=4pq<1$.
When $z\to +\infty$ the probability $P(z)$ decays exponentially, and, more precisely,
$$
P(z) \sim \frac{s^{z}}{\sqrt{\pi (1-s) z}} \ .
$$
\end{cor}

One can obtain higher order asymptotics in the classical way, or by using 
equivalent combinatorical methods as in \cite{GZ2018}

We can be more precise by looking at the time it takes to the honest network to mine $z$ blocks. A longer duration 
than the average $\tau_1=z\tau_0$ leaves extra time for the attacker to build his replacement blockchain, and with this conditional knowledge 
the probability changes. If we define $\kappa = \frac{\tau_1}{z\tau_0}$, we can compute this probability 
$P(z,\kappa)$ and we can also obtain an exact formula using the regularized incomplete Gamma function
$$
Q(s,x)=\frac{\Gamma(s,x)}{\Gamma(x)}
$$
where
$$
\Gamma(s,x)= \int_x^{+\infty} t^{s-1} e^{-t} \, dt
$$
is the incomplete Gamma function.

\begin{theorem} We have
$$
 P(z, \kappa)=1-Q(z, \kappa\, z\, q/p) + \left (\frac{q}{p}\right )^{z} e^{\kappa\, z\,\frac{p-q}{p}} Q(z, \kappa\, z) \ ,
$$
\end{theorem}

\begin{figure}
\centering
\includegraphics[width=1\linewidth]{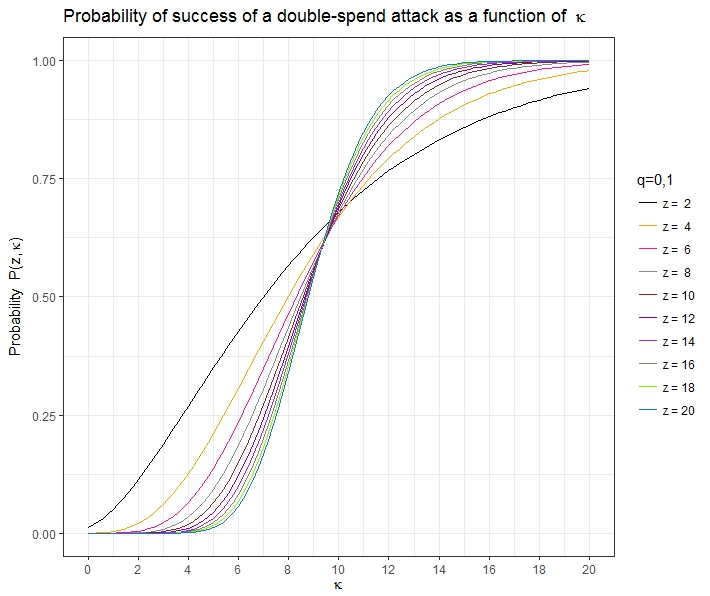}

\smallskip

\caption{Probabilities $P(z, \kappa)$ for $q=0.1$ and distinct values of $z$.}\label{fig:probabilites}

\end{figure}

Figure \ref{fig:probabilites} shows the graphs of $\kappa\mapsto P(z,\kappa)$ for different values of $z$.

\section{Mining profitability.} \label{sec:profitability}
\smallskip

After the study of the security of the protocol, the next important problem is its stability. For a decentralized protocol it is 
fundamental that the interests of the individuals are properly aligned with the protocol rules. In particular, the maximal gain of 
miners should be achieved when following the protocol rules. This is far from obvious, and we know from the 
study of unstable Dynamical Systems that this is hard to achieve. It is somewhat surprising that this has been empirically verified 
since Bitcoin's inception. 

For example, it is by no means obvious that is in the best interest of a miner to publish immediately 
a block that he has validated. He can 
keep it secret and secretly push his advantage, but then he runs the risk that another 
miner publishes a validated block and the 
public blockchain adopts it thus loosing his reward. This type of scenario 
has been discussed since 2012 in bitcointalk forum, created by Nakamoto in 2010.

To answer this question, we first need  to develop a proper profitability model.
As in any business, mining profitability is accounted by the ``Profit and Loss'' per unit of time. The profits of a miner come from 
the block rewards that include the coinbase reward in new bitcoins created, and the transaction fees of the transactions in the block.
The profitability at instant $t>0$ is given by 
$$
\mathbf {PL}(t) =\frac{\mathbf R(t)-\mathbf C(t)}{t}
$$
where $\mathbf R(t)$ and $\mathbf C(t)$ represent, respectively, the rewards and the cost of the mining operation up to time $t$.
If we don't consider transaction fees we have 
$$
\mathbf R(t) = \mathbf N(t) \,  b
$$ 
where $b>0$ is the coinbase reward. If we include transaction fees, the last equation remains true taking the average reward using the classical 
Wald Theorem. 

The random variable $\mathbf C(t)$ representing the cost of mining operations is far more complex to determine since it 
depends on external factors (as electricity costs, mining hardware costs, geographic location, currency exchange rate, etc). 
But, fortunately, we don't need it in the comparison the profitability of different mining strategies as we explain next.

The mining activity is a repetitive and the miners returns to the same initial state after some time, for instant, start mining
a fresh block. A mining strategy is composed by cycles where the miner invariably returns to the initial state. It is a ``game with 
repetition'' similar to those employed by profit gamblers in casino games (when they can spot a weakness that makes the game profitable).
For example, an honest miner starts a cycle each time the network, he or someone else, has validated a new block. 

We denote by $\boldsymbol{\tau}$ the duration of the cycle, and we are interested in \textit{integrable} strategies for which  
$\mathbb E [\boldsymbol{\tau}] <+\infty$ (this means that the cycles 
almost surely end up in finite time). Then it is easy to check,
using the law of large numbers and Wald Theorem, that the long term profitability is given a.s. by the limit 
$$
\mathbf {PL}_\infty =\lim_{t\to +\infty}\frac{\mathbf R(t)-\mathbf C(t)}{t} =\frac{\mathbb E [\mathbf R]-\mathbb E [\mathbf C]}{\mathbb E [\boldsymbol{\tau}]}
$$
As observed before the second cost term is hard to compute, but the revenue term, that we call \textit{revenue ratio}, is in general computable 
$$
\mathbf{\Gamma} = \frac{\mathbb E [\mathbf R]}{\mathbb E [\boldsymbol{\tau}]}
$$
For example, for an honest miner we have  $\mathbb E[\mathbf R] = p.0+q . b=qb$ and  $\mathbb E[\boldsymbol{\tau}] =\tau_0$, and therefore
$$
\mathbf{\Gamma}_H =\frac{qb}{\tau_0} 
$$

We have the fundamental Theorem on comparison of mining strategies 
with the same cost ratio. This is the case when both strategies use the full mining power at all time.

\begin{theorem}[\cite{GPM3}, 2018]\label{thm_comparaison}
We consider two mining strategies $\xi$ and $\eta$ with the same cost by unit of time.
Then  $\xi$ is more profitable than $\eta$ if and only if
$$
\mathbf{\Gamma}_\eta \leq \mathbf{\Gamma}_\xi
$$
\end{theorem}

\section{Protocol stability.}

We can now mathematically study the protocol stability. The following remarkable result (remarkable because it is hard to imagine how Nakamoto could have foreseen it)
validates the proper adjustment of the protocol:

\begin{theorem}[\cite{GPM3}, 2018]\label{thm_stabilite}
In absence of difficulty adjustment, the optimal mining strategy is to publish immediately all mined blocks as soon as they are discovered.
\end{theorem}

We remind that the difficulty of mining adjusts in about every two weeks, so at the same time we spot a weakness of the protocol that we discuss below.

This Theorem holds true for any hashrate of the miner and without any assumption of the type of miners present in the network. It changes nothing that eventually 
there are some dishonest miners in the network.

The proof is simple and is a good example of the power of martingale techniques. For a constant difficulty, the average speed of block 
discovery remains constant and the counting process $\mathbf N(t)$ is a Poisson process with intensity $\alpha=\frac{p}{\tau_0}$ where 
$p$ is the relative hashrate of the miner. The cycle duration $\boldsymbol{\tau}$ is a stopping time and the revenue per cycle equals to 
$\mathbf R=\mathbf N\left(\boldsymbol{\tau}\right)$. Its mean value is then obtained using Doob's stopping time to the martingale 
$\mathbf N(t)-\alpha t$. Finally we get $\mathbf{\Gamma} \leq \mathbf{\Gamma}_H$.

But the Bitcoin protocol does have a difficulty adjustment algorithm that is necessary, in particular during the development phase. Theorem \ref{thm_stabilite}
shows that this is  the only vector of attack. This difficulty adjustment provides a steady monetary creation and ensures that the interblock validation time 
stays at around $10$ minutes. A minimum delay is necessary to allow a good synchronization of all network nodes. If the hashrate accelerates without 
a difficulty adjustment, then the nodes will desynchronise, and many competing blockchains will appear leaving a chaotic state.

\section{Profitability of rogue strategies.}

In view of Theorems \ref{thm_comparaison} and \ref{thm_stabilite}, and in order to decide if a mining strategy is more profitable than 
the honest strategy, we only need to compute the revenue ratio $\mathbf{\Gamma}$ with the difficulty adjustment integrated.
Selfish mining (SM strategy 1) is an example of rogue strategy. Instead of publishing a new block, the miner keeps the block secret and 
tries to build a longer blockchain increasing its advantage. When he makes it public, he will orphan the last mined honest blocks 
and will reap the rewards. To be precise, the attack cycles are defined as follows: the miner starts mining a new block on top of the 
official blockchain. If an honest miner finds a block first, then the cycle ends and he starts over. Otherwise, when he is the first to find a block, 
he keeps mining on top of it and keeping it secret. If before he mines a second block the honest network mines one public block, then he publishes his block 
immediately, thus trying to get a maximum proportion $0<\gamma<1$ of honest miners adopting his block. The propagation is not instantaneous and the efficiency depends 
on the new parameter $\gamma$ which represents his good connectivity to the network. A competition follows, and if the next block is mined on top of the honest 
block, then the selfish miner looses the rewards of this block and the attack cycle ends. If he, or his allied honest miners, mine the next block, then they publish it
and the attack cycle ends again. If the attacker succeeds in mining two consecutive secret blocks at the beginning, then he continues 
working on his secret blockchain until he has only one block of advantage with respect to the public blockchain. In that case, he doesn't run any 
risk of being joined by the public blockchain and publishes all his secret blockchain, thus  reaping all the rewards and ending the attack cycle again. In few words,
the rogue miner spends most of his time replacing honest blocks by those that he mined in advance and kept secret. The mean duration 
$\mathbb E [\boldsymbol{\tau}]$ of the attack cycle is obtained as a variation of the following result about Poisson processes.

\begin{prop}[Poisson races]
Let $\mathbf N$ and $\mathbf N'$ be two independent  Poisson processes with respective parameters $\alpha$ and 
$\alpha'$, with $\alpha' < \alpha$ and $\mathbf N(0) = \mathbf N'(0) = 0$. Then, the stopping time 
\begin{equation*}
\boldsymbol{\sigma} = \inf \{ t>0 ; \mathbf N(t) = \mathbf N'(t) + 1\}
\end{equation*} 
is almost surely finite, and we have
\begin{equation*}
\mathbb E [\boldsymbol{\sigma}] = \frac{1}{\alpha - \alpha'}\, ,  \ \ \mathbb E[\mathbf N'(\boldsymbol{\sigma})] =  \frac{\alpha'}{\alpha - \alpha'} \, , 
\ \ \mathbb E[\mathbf N(\boldsymbol{\sigma})] =  \frac{\alpha}{\alpha - \alpha'} \ . 
\end{equation*}
\end{prop}
The proof is a simple application of Doob's Stopping Time Theorem.
Here, $\mathbf N$ and $\mathbf N'$ are the counting processes of blocks mined by the honest miners and the attacker. 
To finish, we must compute the intensities $\alpha$ and $\alpha'$. At the beginning we have $\alpha =\alpha_0 = \frac{p}{\tau_0}$
and $\alpha' =\alpha'_0 = \frac{q}{\tau_0}$ where $p$ is the apparent hashrate of the honest miners and $q$ the one of the attacker. 
But the existence of a selfish miner perturbs the network and slows down the production of blocks. Instead of having one block 
for each period $\tau_0$, the progression of the official blockchain is of $\mathbb E [N(\boldsymbol{\tau})\vee N'(\boldsymbol{\tau})]$ 
blocks during $\mathbb E [\boldsymbol{\tau}]$. After validation of $2016$ official blocks, this triggers a difficulty adjustment that 
can be important. The new difficulty is obtained from the old one by multiplication by a factor $\delta <1$ given by 
$$
\delta = 
\frac{\mathbb E [N(\boldsymbol{\tau})\vee N'(\boldsymbol{\tau})]\, \tau_0}{\mathbb E [\boldsymbol{\tau}]}
$$
After the difficulty adjustment, the new mining parameters are $\alpha =\alpha_1 = \frac{\alpha_0}{\delta}$
and $\alpha' =\alpha'_1 = \frac{\alpha'_0}{\delta}$. The stopping time $\boldsymbol{\tau}$ and the parameter $\delta$ 
can be computed using the relation $|N(\boldsymbol{\tau}) - N'(\boldsymbol{\tau})|=1$. This can be used to compute the 
revenue ratio of the strategy \cite{GPM3}. This analysis can also checked by mining simulators.

An alternative procedure consists in modelling the network by a Markov chain where the different states correspond to different degree 
of progress of the selfish miner. Each transition corresponds to a revenue increase $\boldsymbol{\pi}$ and $\boldsymbol{\pi}'$ for the 
honest and selfish miner. By another application of the Law of Large Numbers we prove that the long term apparent hashrate of the 
strategy, defined as the proportion of mined blocks by the selfish miner compared to the total number of blocks, is given by 
the formula
$$
q' = 
\frac{\mathbb E [\boldsymbol{\pi}']}
{\mathbb E [\boldsymbol{\pi}]+\mathbb E [\boldsymbol{\pi}']}
$$

The expectation is taken relative to the stationnary probability that exists 
because the Markov chain is transitive and recurrent. Indeed, the Markov chain is essentially a random walk on $\mathbb{N}$
partially reflexive on $0$. The computation of this stationnary probability proves the following Theorem:

\begin{theorem}[\cite{ES14}, 2014]\label{enc}
The apparent hashrate of the selfish miner is
$$
q' = \frac{((1 + pq) (p - q) + pq) q - (1 - \gamma) p^2 q (p - q)}{p^2 q + p - q}
$$
\end{theorem}

The results from \cite{GPM3} and \cite{ES14} obtained by these different methods, are compatible. The revenue ratio $\boldsymbol{\Gamma}_{1}$
and the apparent hashrate $q'$ are related by the following equation
$$
\boldsymbol{\Gamma}_{1} = q'\frac{b}{\tau_0}
$$
But the first analysis is finer since it does explain the change of  profitability regime after the difficulty adjustment. In particular, it allows to 
compute the duration before running into profitability for the attacker. The selfish miner starts by losing money, then after the difficulty adjustment 
that favors him, starts making profits. For example, with $q=0.1$ and $\gamma=0.9$, he needs to wait $10$ weeks in order to be profitable. This 
partly explains why such an attack has never been observed in the Bitcoin network.

Theorem \ref{thm_comparaison} gives an explicit semi-algebraic condition on the parameters, namely  $q'> q$, that determines the values of the parameters 
$q$ and $\gamma$ for which the selfish mining strategy is more profitable than honest mining.

Theorem \ref{thm_stabilite} shows that the achilles' heel of the protocole is the difficulty adjustment formula. This formula is supposed to 
contain the information about the total hashrate, but in reality it ignores the orphan blocks. The authors proposed a solution that incorporates this 
count, and this solves the stability problem of the protocol \cite{GPM3}.

There are other possible block-withholding strategies that are variations of the above strategy \cite{NKMS2016}. These are 
more agressive strategies. In the initial situation where the attacker succeeds to be two blocks ahead, instead of publishing the whole 
secret chain when he is only one block ahead, he can wait to be caught-up to release his blocks and then starts a final competition between the two 
competing chains. The attack cycle ends when the outcome is decided. This is the ``Lead Stubborn Mining'' (LSM, strategy $2$). In this strategy 
it is important that the miner regularly publishes his secret blocks with the same height of the official blockchain, to attract part of the honest 
miners in order to take out hashrate from the pool of honest miners. Also in this way, even if he looses the final competition he will succeed 
in incorporating some of his blocks in the official blockchain and reap the corresponding rewards. 

Another even more agressive variation consists in waiting not to be 
caught up but 
to be behind one block. This is the "Equal Fork Stubborn Mining Strategy" (EFSM, strategy $3$). 
Here again, it is important to publish secret blocks regularly. Finally, the authors have 
considered another more agressive variation where the stubborn miner follows EFSM but then 
doesn't stop when he is one block behind. He keeps on mining until his delay becomes greater 
than a threshold $A$ or until he successfully comes from behind, catches-up and finally 
takes the advantage over the official blockchain. 

This strategy seems desperate, because the official blockchain progress is faster, on average. But in case of catching-up the selfish 
miner wins the jackpot of all the blocks he replaces. This is the ``A-Trailing Mining'' strategy (A-TM, strategy $4$). The authors of \cite{NKMS2016}
conduct a numerical study of profitability by running a Montecarlo simulation and compare the profitability of the different strategies for 
different parameter values $(q,\gamma)$. But we can find closed form formulas for the revenue ratio of all these strategies using the precedent 
martingale approach.

\begin{theorem}(\cite{GPM3,GPM4,GPM5,GPM6})
We have
\begin{align*}
\frac{\mathbf \Gamma_{1}}{\mathbf \Gamma_H}&=\frac{(1 + pq) (p - q) + pq - (1 - \gamma) p^2  (p - q)}{p^2 q + p - q} \\
\frac{\mathbf\Gamma_{2}}{\mathbf\Gamma_H}&=   \frac{p + pq - q^2}{p + pq -q} - \frac{p (p - q)  f(\gamma, p,q)}{p + pq - q} \\
\frac{\mathbf\Gamma_{3}}{\mathbf\Gamma_H}&=  \frac{1}{p} - \frac{p-q}{pq}  f(\gamma, p, q) \\
\frac{\mathbf \Gamma_{4}}{\mathbf \Gamma_H}&= \frac{1 + \frac{(1 - \gamma) p (p - q)}{(p + pq - q^2) [A + 1]}  \left(
  \left( [A - 1] + \frac{1}{p}  \frac{P_A (\lambda)}{[A + 1]} \right) \lambda^2
  - \frac{2}{\sqrt{1 - 4 (1 - \gamma) pq} + p - q}  \right)}
  {\frac{p + pq - q}{p + pq -q^2} + \frac{(1 -
  \gamma) pq}{p + pq - q^2}  (A + \lambda)  \left( \frac{1}{[A + 1]} - \frac{1}{A + \lambda} \right)}
\end{align*}
with 
$$
f(\gamma, p,q)=\frac{1 - \gamma}{\gamma} \cdot \left (1 - \frac{1}{2q} (1-\sqrt{1-4(1 - \gamma) pq})\right )
$$
and $\lambda = q/p$, $[n]=\frac{1-\lambda^n}{1-\lambda}$ pour  $n \in \mathbb{N}$, 
$P_A (\lambda) = \frac{1 - A \lambda^{A - 1} + A \lambda^{A + 1} - \lambda^{2 A}}{(1 - \lambda)^3}$
\end{theorem}

\begin{figure}
\centering
\includegraphics[width=0.8\linewidth]{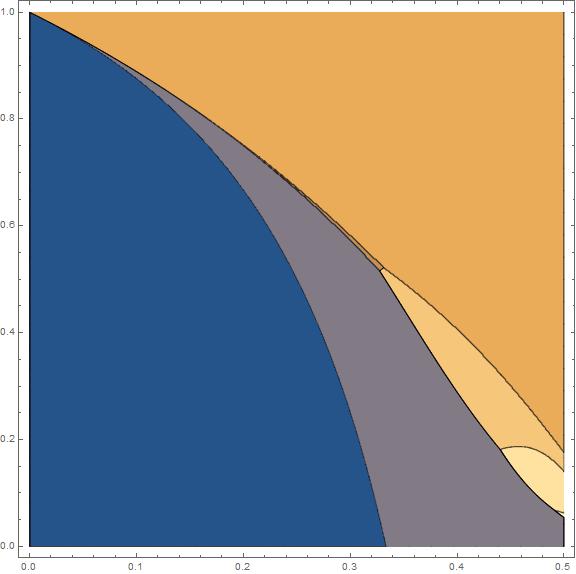}

\smallskip

\caption{Comparison of HM, SM, LSM, EFSM and A-TSM.}\label{fig:dominance}

\end{figure}

We can plot the parameter regions where each startegy is the best one (Figure \ref{fig:dominance}).
The Catalan numbers appear naturally in the computations.
$$
 C_n = \frac{1}{2n+1} \binom{2n}{n} = \frac{(2 n) !}{n! (n + 1) !}  \ .
$$ 
For this reason their generating function appears in the formulas
$$
 C (x) = \sum_{n=0}^{+\infty} C_n x^n = \frac{1-\sqrt{1-4x}}{2x}
$$
We observe that $\sqrt{1-4pq}=p-q$ and $C(pq)=1/p$, and this justifies the definition of new probability distributions that arise in the proofs.
\begin{definition}
A discrete random variable  $X$ taking integer values follows a Catalan distribution of the first type if we have, for $n\geq 0$,
 $$
 \mathbb P [X=n] = C_{n} p(pq)^n  \ .
 $$
It follows a Catalan distribution of the second type if  $\mathbb P [X=0]=p$ and for $n\geq 1$,
 $$
 \mathbb P [X=n]= C_{n-1} (pq)^n  \ .
 $$
It follows a Catalan distribution of the third type if $\mathbb P [X=0]=p$, $\mathbb P [X=1]=p q+p q^{2}$ 
and for $n\geq 2$,
 $$
 \mathbb P [X=n]= p q^{2} C_{n-1} (pq)^{n-1}  \ .
 $$
\end{definition}

\section{Dyck words.}

We can recover this results by a direct combinatorical approach representing each attack cycle by a Dyck word.

\begin{definition}
A Dyck word is a word built from the two letter alphabet $\{S,H\}$ which contains as many S letters as H letters, and such that any prefix word contains 
more or equal S letters than H letters. We denote $\mathcal D$ the set of Dyck words, and for  $n\geq 0$, $\mathcal D_n$ the subset of Dyck worlds of 
length $2n$.
\end{definition}

The relation to Catalan numbers is classical: the cardinal of 
$\mathcal D_n$ is $C_n$. We can encode attack cycles by a chronologic succession of block discoveries (disregarding if the blocks are made public or not). 
For a selfish block we use the letter S (for ``selfish'') and for the honest blocks the letter H (for ``honest'').

The link between the selfish mining strategy and Dyck words is given by the following proposition:

\begin{prop}
The attack cycles of the SM strategy are H, SHH, SHS, and SSwH where $w\in \mathcal D$. 
\end{prop}

At the end of the cycle, we can summarise and count the total number of official blocks, say $L$, and how 
many of these blocks were mined by the attacker, say $Z$. Then, for strategy 1 (SM), the random variable $L-1$ follows a Catalan distribution 
of the third type, and except for some particular cases (when $L<3$), we always have $L=Z$. The apparent hashrate $q'$ is then given by the formula:
$$
q' = \frac{\mathbb E [Z]}{\mathbb E [L]}
$$
We can then directly recover Theorem \ref{enc} by this simpler combinatorical procedure \cite{GPM6}.
The other rogue strategies can be studied in a similar way. The Catalan distribution of the first type arises in the 
study of the strategy EFSM (strategy $3$), and the one of the second type for the strategy LSM (strategy $2$).
We can then recover all the results given by the Markov chain analysis. Unfortunately we cannot recover the 
more finer results obtained by martingales techniques.

This sort analysis applies to other Proof-of-Work cryptocurrencies, and to Ethereum that has a more complex reward system and a different difficulty adjustment formula \cite{GPM7}.

\section{Nakamoto double spend revisited.}

We come back to the fundamental double spend problem from Nakamoto Bitcoin paper 
discussed in section \ref{sec:double_spend}. In that section, 
we computed the probability of success of a double spend. But now, with the profitability 
model knowledge from section \ref{sec:profitability}, we can study its profitability and get better 
estimates on the number of confirmations that are safe to consider a paiement definitive.
The double spend strategy as presented in \cite{Nakamoto08} is unsound because there is a non-zero probability of failure
and in that case, if we keep mining in the hope of catching-up from far behind the official blockchain, we have a positive 
probability of total ruin. Also the strategy is not integrable since the expected duration of the attack is infinite.
Thus, we must obviously put a threshold to the unfavorable situation where we are lagging far behind the official blockchain. 

We assume that the number of confirmations requested by the recipient of the transaction is $z$ and we assume that we are never 
behind $A\geq z$ blocks of the official blockchain. This defines an integrable strategy, The $A$-Nakamoto double spend strategy. Putting aside technical details about premining, the probability of success of this strategy is a modification of the probability
from Theorem \ref{thm:probability}

\begin{theorem} [\cite{GPM8}, 2019]
After $z$ confirmations, the probability of success of an  $A$-Nakamoto double 
spend  is
 $$
 P_A(z)=\frac{P(z)-\lambda^A}{1-\lambda^A} \ 
 $$
where $P(z)$
is the probability from Theorem \ref{thm:probability} and $\lambda=q/p$.
\end{theorem}

If  $v$ is the amount to double spend, then we can compute the revenue ratio 
$\mathbf{\Gamma}_A = \mathbb E [\mathbf R]/\mathbb E [\boldsymbol{\tau}]$.

\begin{theorem}[\cite{GPM8}, 2019]
With the previous notations, the expected revenue and the expected duration 
of  the $A$-Nakamoto double spend strategy is
    \begin{align*}
	\mathbb E[\mathbf{R}_A]/b &=\frac{q z}{2 p} I_{4pq}(z,1/2)
	- \frac{A\lambda^{A}}
	{p(1-\lambda)^3 [A]^2} I_{(p-q)^2}(1/2,z) \\
	&\ \ \ +\frac{2-\lambda+\lambda^{A+1}}{(1-\lambda)^2[A]}
	\frac{p^{z-1}q^{z}}{B(z,z)} 
	+ P_A(z) \left (\frac{v}{b}+1\right ) \\
    \mathbb E[\mathbf{T}_{A}]/\tau_0 &=
	\frac{z}{2p} I_{4pq}(z,1/2) +
	\frac{A}{p(1-\lambda)^2[A]}
	I_{(p-q)^2}(1/2,z)  \\
    &\ \ \ - \frac{p^{z-1} q^z}{p(1-\lambda) \, B(z,z)} 
    +\frac{1}{q}
    \end{align*}  
with the notation $[n]=\frac{1-\lambda^n}{1-\lambda}$ for an integer $n\geq 0$, and $B$ 
is the classical Beta function.
\end{theorem}
 
\begin{figure}
\centering
\includegraphics[width=0.8\linewidth]{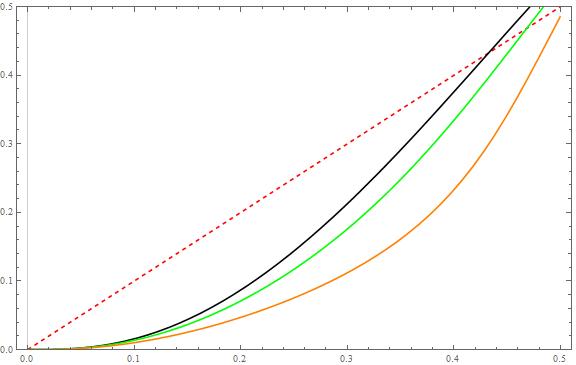}

\smallskip

\caption{Graph of $\Gamma_A$, for $z = 2$,  $v=b$ and  $A = 3,5,10$.}

\end{figure}

In principle a powerful miner does not have an interest in participating in a large double spend, since
doing so will undermine the foundations of his business. For a small miner with relative hashrate $0<q<<1$ 
we can estimate the minimal amount of a double spend to be profitable. For this we only need to use 
the inequality from Theorem \ref{thm_comparaison}, $\mathbf{\Gamma}_A  \geq \mathbf{\Gamma}_H=qb/\tau_0$, and take the asymptotics 
$q\to 0$ (with $A$ and $z$ being fixed, but the final result turns out to be independent of $A$). 

\begin{corollary}
When $q\to 0$ the minimal amount $v$ for an Nakamoto double spend with $z\geq 1$ confirmations
is 
$$
v\geq \frac{q^{-z}}{2\, \binom{2z-1}{z}} \, b =v_0 \ .
$$
\end{corollary}

For example, in practice, with a $10$\% hashrate, $q=0.01$, and only one confirmation, $z=1$, 
we need to double spend more than $v_0/b=50$ coinbases. 
With the actual coinbase reward of $b=12.5 \, \bitcoinA$ bitcoins and the actual prize over $8.300$ euros, this represents more than $5$  millions euros. 

Hence for all practical purposes and normal 
amount transactions, only one confirmation is enough to consider the transaction definitive.

\medskip

\textbf{Conclusions.}
Bitcoin provides a good example of the universality of Mathematical applications and its potential to impact our society. With the glimpse we have given, we hope to have convinced our colleagues that the Bitcoin protocol also motivates some exciting Mathematics.

\bibliographystyle{plain}
\bibliography{bitcoin}

\end{document}